\magnification\magstep1
\baselineskip = 18pt
\def\n{\noindent}
\magnification\magstep1

\baselineskip = 18pt
\def\n{\noindent}

\def\pf{\medskip{\noindent{\bf Proof. }}}
\def \rat{ {\rm Q}\kern-.65em {}^{{}_/ }}

 \baselineskip = 18pt
 \def\n{\noindent}
 \overfullrule = 0pt
 \def\qed{{\hfill{\vrule height7pt width7pt
depth0pt}\par\bigskip}} 
 
 \magnification\magstep1

\def \rat{ {\rm Q}\kern-.65em {}^{{}_/ }}

\overfullrule = 0pt
\def\pf{\medskip{\noindent{\bf Proof. }}}

\def\vp{\varepsilon}

\overfullrule = 0pt
\def\pf{\medskip{\noindent{\bf Proof. }}}

\centerline{\bf Exact operator spaces}
\bigskip
\centerline{by}\bigskip
\centerline{Gilles Pisier\footnote*{Supported in part by the NSF}}
\centerline{Texas A\&M University}
\centerline{and}
\centerline{Universit\'e Paris VI}
\n Plan

\n Introduction

\n \S 1. Exact operator spaces.

\n \S 2. Ultraproducts.

\n \S 3. How large can $d_{SK}(E)$ be?

\n \S 4. On the set of $n$-dimensional operator spaces.

\n \S 5. On the dimension of the containing matrix
space.\vfill\eject

\n {\bf Introduction}

In this paper, we study {\it operator spaces\/} in the sense of the theory
developed recently by Blecher-Paulsen [BP] and Effros-Ruan [ER1]. By an
operator space, we mean a closed subspace $E\subset B(H)$, with $H$
Hilbert. In the category of operator spaces, the morphisms are the {\it
completely bounded maps\/} for which we refer the reader to [Pa1]. Let
$E\subset B(H)$, $F\subset B(K)$ be operator spaces ($H,K$ Hilbert). A map
$u\colon \ E\to F$ is called completely bounded ($c.b.$ in short) if
$$\sup_{n\ge 1}\|I_{M_n}\otimes u\|_{M_n(E)\to M_n(F)}<\infty$$
where $M_n(E)$ and $M_n(F)$ are equipped with the norms induced by
$B(\ell^n_2(H)$) and $B(\ell^n_2(K))$ respectively. We denote
$$\|u\|_{cb}=\sup_{n\ge 1}\|I_{M_n}\otimes u\|_{M_n(E)\to M_n(F)}.$$
The map $u$ is called a complete isomorphism if it is an isomorphism and if
$u$ and $u^{-1}$ are $c.b.$. We say 
that $u\colon\ E\to F$ is a complete
isometry if for each $n\ge 1$ the map $I_{M_n}\otimes u\colon \ M_n(E)\to
M_n(F)$ is an isometry. We refer to [Ru, ER2-7, B1, B2] 
for more information on the rapidly developing Theory of Operator
Spaces.

We will be mainly concerned here with the ``geometry'' of {\it finite
dimensional\/} operator spaces. In the Banach space category, it is
well known that every separable space embeds isometrically into
$\ell_\infty$. Moreover, if $E$ is a finite dimensional normed space then
for each $\vp>0$, there is an integer $n$ and a subspace $F\subset
\ell^n_\infty$ which is $(1+\vp)$-isomorphic to $E$, i.e. there is an
isomorphism $u\colon \ E\to F$ such that $\|u\|\ \|u^{-1}\|\le 1+\vp$. Here
of course, $n$ depends on $\vp$, say $n=n(\vp)$ and usually (for instance
if $E = \ell^k_2$) we have $n(\vp)\to \infty$ when $\vp\to 0$.

Quite interestingly, it turns out that this fact is not valid in the
category of operator spaces:\ although every operator space embeds
completely isometrically into $B(H)$ (the non-commutative analogue of
$\ell_\infty$) it is not true that a finite dimensional operator space must
be close to a subspace of $M_n$ (the non-commutative analogue of
$\ell^n_\infty$) for some $n$. The main object of this paper is to study
this phenomenon.

 We will see that this phenomenon is very closely related to the remarkable
work of E.~Kirchberg on exact $C^*$-algebras. We will show that some of
Kirchberg's ideas can be  developed in a purely ``operator space'' setting.
Our main result in the first section  is Theorem 1, 
which can be stated as follows.

\n
Let $B=B(\ell_2)$ and let $K\subset B$ be the ideal of
all the compact operators on $\ell_2$. 

\n If $X,Y$  are
operator spaces, we denote by $X\otimes_{\min} Y$
their minimal (or spatial) tensor product. If $X\subset B(H)$ and
$Y\subset B(K)$, this is just the completion of the linear tensor product
 $X\otimes Y$ for the norm induced by $B(H\otimes_2 K)$.

\n Let $\lambda\ge 1$ be a fixed
constant.

\n The following properties of an
operator space $X$ are equivalent:

(i) The sequence 
$$\{0\} \to K\otimes_{\rm min}X \to B\otimes_{\rm min}X \to (B/K)
\otimes_{\rm min} X\to \{0\}$$
is exact  and the map  
$$T_X\colon \ (B\otimes_{\rm min} X)/(K\otimes_{\rm min} X) \to (B/K)
\otimes_{\rm min} X $$
has an inverse $T_X^{-1}$ with norm
$\|T_X^{-1}\| \le \lambda$. 

(ii) for each $\epsilon >0$ and each finite dimensional subspace
$E\subset X$, there is an integer $n$ and a subspace $F \subset M_n$ such that
$d_{cb}(E,F) <\lambda +\epsilon$.

\n Here $d_{cb}(E,F)$ denotes the $c.b.$ analogue of the Banach-Mazur
distance (see (0) below for a precise definition.) We will denote
by $d_{SK}(E)$ the infimum of $d_{cb}(E,F)$ when $F$ 
runs over all
operator spaces $F$ which are subspaces of $M_k$ for some
integer  $k$.

 One of the main results in section 2 can be stated as follows (see Theorem~7
below). 

\n Consider $F\subset M_k$ with $\dim F=n$ and $k\ge n$ arbitrary, then for any linear
isomorphism $u\colon \ \ell^n_\infty\to F^*$ we have
$$\|u\|_{cb}\|u^{-1}\|_{cb} \ge n[2(n-1)^{1/2}]^{-1}.$$
In particular this is $>1$ for any $n\ge 3$. Here the space $F^*$ is the
dual of $F$ with its ``dual operator space structure'' as explained in
[BP, ER1, B1, B2].

\n Equivalently, if we denote
by $E_1^n$ the  operator space dual of $\ell_\infty^n$,
(this is denoted by $\max(\ell_1^n)$ in [BP]) then we have
$$d_{SK}(E_1^n)\ge {n\over 2\sqrt{n-1}}.$$
We also show a similar estimate for the space 
which is denoted by $R_n+C_n$ in [P1]. Moreover, we show that
the $n$-dimensional operator Hilbert space $OH_n$ (see [P1])
satisfies
$$d_{SK}(OH_n)\ge ({n\over 2\sqrt{n-1}})^{1/2}.$$
These estimates are asymptotically sharp in the sense that
$d_{SK}(E_1^n)$ and $d_{SK}(R_n+C_n)$ are $O(n^{1/2})$ and 
$d_{SK}(OH_n)$ is $O(n^{1/4})$ when $n$ goes to infinity.
 
Later on in the paper, we show that the operator space analogue of the
``Banach Mazur compactum'' is not compact and we prove various estimates
related to that phenomenon. (The noncompactness itself was known, at least
to Kirchberg.) We will include
several simple facts on ultraproducts of finite dimensional operator spaces
which are closely connected to the discussion of ``exact'' operator spaces
presented in section~1. Let us denote by
$OS_n$ the set of all $n$-dimensional operator spaces. We consider that two
spaces $E,F$ in $OS_n$ are the same if they are completely isometric. Then
the space $OS_n$ is a metric space when equipped with the distance
$$\delta(E,F) = \log d_{cb}(E,F).$$
We include a proof that $OS_n$ is complete but not compact (at least if
$n\ge 3$) and we give various related estimates. 
As pointed out to me by Kirchberg, it seems to be an open
problem whether $OS_n$ is a {\it separable\/} metric space.

In passing, we recall that in [P1] we proved that $d_{cb}(E, OH_n) \le
n^{1/2}$ for any $E$ in $OS_n$ and therefore that
$$\sup\{d_{cb}(E,F)\mid E,F\in OS_n\}=n.$$
Actually, that supremum is attained on the subset 
$HOS_n\subset OS_n$ formed of all the 
{\it Hilbertian\/} operator spaces (i.e. those which, as normed spaces, are
isometric to the Euclidean space $\ell^n_2$). 
We also show that (at least
for $n\ge 3$) $HOS_n$ is a closed but non compact  
subset of $OS_n$.
 Perhaps the subset  $HOS_n$  is not even
separable. 

In section 5, we show the following result. Let $E$ be any operator space and let $C\ge 1$ be a
constant. Fix an integer $k\ge 1$. Then there is a compact set $T$ and a subspace
$F\subset C(T)\otimes_{\rm min} M_k$ such that
$d_{cb}(E,F) \le C$ iff  for any operator space
$X$ and any $u\colon \ X\to E$ we have $$\|u\|_{cb} \le
C\|u\|_k,$$ where $\|u\|_k=\|u\|_{M_k(X)\to M_k(E)}$.

\n {\bf Notation:} Let $(E_m)$ be a sequence of operator spaces.
We denote by $\ell_{\infty}\{E_m\}$ the direct sum in the sense
of $\ell_{\infty}$ of the family  $(E_m)$. 
As a Banach space,
this means
 that 
$\ell_{\infty}\{E_m\}$ is the set of all  sequences
$x=(x_m)$ with $x_m\in E_m$ for all $m$ with  
$\sup_m \|x_m\|_{E_m} <\infty$ equipped with the norm
$\|x\|=\sup_m \|x_m\|_{E_m} $. The operator space structure on
$\ell_{\infty}\{E_m\}$ is defined by the identity
$$\forall n \quad M_n(\ell_{\infty}\{E_m\})=\ell_{\infty}\{M_n(E_m)\}
.$$ Equivalently, if $E_m \subset B(H_m)$ (completely isometrically)
then $\ell_{\infty}\{E_m\}$ embeds (completely isometrically)
into $B(\oplus_m H_m)$ as block diagonal operators.

\n We will use several times the observation that
if $F$ is an other operator space then
$\ell_{\infty}\{E_m\}\otimes_{\min} F$ embeds completely isometrically
in the natural way into $\ell_{\infty}\{E_m\otimes_{\min} F \}$. In particular,
if $F$ is finite dimensional these spaces can be completely
isometrically identified.

\n {\bf Acknowlegement.} I am very grateful to E. Kirchberg
 for copies of his papers and for several very 
stimulating conversations.

\vfill\eject

\n {\bf \S 1. Exact operator spaces.}

Let $E,F$ be operator spaces. We denote
$$d_{cb}(E,F) = \inf\{\|u\|_{cb}\|u^{-1}\|_{cb}\}\leqno(0)$$
where the infimum runs over all isomorphisms $u\colon \ E\to F$. If $E,F$
are not completely isomorphic we set $d_{cb}(E,F)=\infty$. This is the
operator space version of the Banach Mazur distance. We will study the
smallest distance of an operator space $E$ to a subspace of the space $K
=K(\ell_2)$ of all compact operators on $\ell_2$. More precisely, this is
defined as follows
$$d_{SK}(E) = \inf\{d_{cb}(E,F)\mid F\subset K\}.\leqno (1)$$
Let $F$ be a finite dimensional subspace of $K$. By an entirely classical
perturbation argument one can check that for each $\vp>0$ there is an
integer $n$ and a subspace $\widetilde F\subset M_n$ such that
$d_{cb}(F,\widetilde F)<1+\vp$. It follows that for any finite dimensional
operator space $E$ we have
$$d_{SK}(E) = \inf\{d_{cb}(E,F)\mid F\subset M_n, \quad n\ge 1\}.\leqno
(1)'$$
In his remarkable work on exact $C^*$-algebras
 (cf. [Ki]) Kirchberg
introduces a quantity which he denotes by ${locfin} (E)$ for any operator
space $E$. His definition uses completely positive unit preserving maps.
The number $d_{SK}(E)$ appears as the natural ``$c.b.$'' analogue of
Kirchberg's ${locfin}(E)$. Note that $d_{SK}(E)$ is clearly an invariant of
the operator space $E$ and we have obviously
$$d_{SK}(E) \le d_{SK}(F) d_{cb}(E,F)\leqno (1)''$$
for all operator spaces $E,F$.

Let $X$ be an operator space. We will say that $X$ is exact if the sequence
$$\{0\} \to K\otimes_{\rm min}X \to B\otimes_{\rm min}X \to (B/K)
\otimes_{\rm min} X\to \{0\}\leqno (2)$$
is exact.

Note in particular that every finite dimensional space is trivially exact.
Following Kirchberg, we will measure the ``degree of exactness'' of $X$ via
the number $ex(X)$ defined as follows:\ we consider the map
$$T_X\colon \ (B\otimes_{\rm min} X)/(K\otimes_{\rm min} X) \to (B/K)
\otimes_{\rm min} X\leqno (3)$$
associated to the exact sequence (2) and we define
$$ex(X) = \|T^{-1}_X\|.\leqno (4)$$
Clearly $ex(X)$ is a (completely) isomorphic invariant of $X$ in the
following sense:\ if $X$ and $Y$ are completely isomorphic operator spaces
we have
$$ex(X) \le ex(Y) d_{cb}(X,Y).\leqno (4)'$$
The main result of this section is the following which is proved by adapting in
a rather natural manner the ideas of Kirchberg [Ki]. One simply
needs to substitute everywhere in his argument ``completely positive
unital'' by ``completely bounded'' and to keep track of the $c.b.$ norms.
The resulting proof is very simple.

\proclaim Theorem 1. For every finite dimensional operator space $E$, one
has
$$ex(E) = d_{SK}(E).\leqno (5)$$
More generally, for any operator space $X$
$$ex(X) = \sup\{d_{SK}(E)\mid E\subset X, \dim E<\infty\},\leqno (6)$$
and $X$ is exact iff the right side of (6) is finite.\medskip

\n {\bf Remarks.} (i)\ If $X$ is a $C^*$-algebra then the maps appearing
in (2) are ($C^*$-algebraic) representations. Recall that a representation
necessarily has closed range and becomes isometric when we pass to the
quotient modulo its kernel (cf. e.g. [Ta,\ p.22]). Hence if $X$ is a $C^*$-algebra (2) is exact
iff the kernel of the map
$$B\otimes_{\rm min} X \to (B/K) \otimes_{\rm min} X$$
coincides with $K\otimes_{\rm min} X$.

\n By a known argument, a sufficient condition for this to hold is a certain
``slice map'' property (cf. [W2]) which is a consequence of the $CBAP$ (see [Kr]).
(Recall that $X$ has the $CBAP$ if the identity on $X$ is a pointwise limit
of a net of finite rank maps $u_i\colon \ X\to X$ with
$\sup_i\|u_i\|_{cb}<\infty$.)

Thus it is known that the reduced $C^*$-algebra of the free group ${\bf
F}_N$ with $N$ generators $(N\ge 2)$ is exact,
 because by [DCH] it has the
$CBAP$. On the other hand, it is known ([W1]) that the full $C^*$-algebras
$C^*({\bf F}_N)$ are not exact. This can also be derived from Theorem~7
below after noticing that the space $E^n_1 = (\ell^n_\infty)^*$ appearing
in
Theorem~7 is completely isometric to a subspace of $C^*({\bf F}_n)$. By the
same argument (using Corollary~10 below) if an operator space $E$ is
completely isometric to the space $OH_n$ introduced in [P1], then the
$C^*$-algebra generated by $E$ is not exact if $n\ge 3$.

\n (ii)\ As explained to me by Kirchberg, if $X$ is a $C^*$-algebra then we
have $ex(X)<\infty$ iff $ex(X) =1$. Indeed since (3) is a representation,
it is isometric if it is injective. This shows if a $C^*$-algebra is
completely isomorphic to an exact operator space then it is exact as a $C^*$-algebra.

We will
use the following simple fact.

\proclaim Lemma 2. Let $E$ be a separable operator space. There are
 operators $P_n\colon \ E\to M_n$ such that\medskip
\item{(i)} $\|P_n\|_{cb} \le 1$ for all $n$.
\item{(ii)} The embedding $J\colon \ E\to \ell_\infty\{M_n\} 
$ defined by $J(x) = (P_n(x))_{n\in {\bf N}}$ is a complete isometry.
\item{(iii)}  For all $k\le n$, there is a map $a_{kn}\colon \ M_n\to M_k$
with $\|a_{kn}\|_{cb}\le 1$ such that $P_k = a_{kn}P_n$.
\item{(iv)} Assume $E$  finite dimensional. Then
  for some $n_0\ge 1$, the maps $P_n$ are injective for all $n\ge
n_0$.\medskip

\n {\bf Proof.} We can assume $E\subset B(\ell_2)$. Then let
$Q_n\colon \ B(\ell_2)\to M_n$ be the usual projection (defined by
$Q_n(e_{ij})=e_{ij}$ if $i,j\le n$ and $Q_n(e_{ij})=0$ otherwise).
Let $P_n = Q_{n|E}$. Then (i), (ii), (iii) and (iv) are immediate.
\qed

The point of the preceding lemma is that we can write for all $N\ge 1$ and
all $(a_{ij})$ in $M_N(E)$
$$\|(a_{ij})\|_{M_N(E)} = \lim_{n\to \infty}\uparrow
\|(P_n(a_{ij}))\|_{M_N(M_n)}.\leqno (7)$$
Indeed, by (ii) we have
$$\|(a_{ij})\|_{M_N(E)}  = \sup_n \|(P_n(a_{ij}))\|_{M_N(M_n)}$$
and by (iii) this supremum is monotone nondecreasing, whence (7).

The following two lemmas are well known to specialists.

\proclaim Lemma 3. If $X,Y$ are exact operator spaces and if $X\subset Y$,
then $ex(X) \le ex(Y)$.

\n {\bf Proof.} We will identify $B\otimes_{\rm min} X$ (resp.
$(B/K)\otimes_{\rm min} X)$ with a subspace of $B\otimes_{\rm min} Y$
(resp. $(B/K)\otimes_{\rm min}Y)$. Consider $u$ in the open unit ball of
$(B/K)\otimes_{\rm min} X$. By definition of $ex(Y)$, there is an element
$v$ in $B\otimes_{\rm min}Y$ such that $\|v\| < ex(Y)$ and if $q\colon \
B\otimes_{\rm min} Y\to (B/K) \otimes_{\rm min} Y$ is the canonical mapping,  we have $q(v) = u$. On
the other hand, since $X$ is exact we know there is a $\tilde u$ in
$B\otimes_{\rm min} X$ such that $q(\tilde u) = u$. Note that by exactness
$\hbox{Ker } q = K\otimes_{\rm min } Y$, hence $v-\tilde u \in
K\otimes_{\rm min}
Y$. Let $p_n$ be an increasing sequence of finite rank projections in $B$
tending to the identity (in the strong operator topology). Consider the mapping $\sigma_n\colon \ B\to B$
defined by $\sigma_n(x) =   (1-p_n)x (1-p_n)$. Clearly
$\|\sigma_n\|_{cb}\le 1$ and for all $x$ in $B$ we have $\sigma_n(x) -x\in
K$. Moreover, for all  $x$ in $K$  we have $\|\sigma_n(x)\|\to 0$. More generally,
 by equicontinuity
$\|\sigma_n\otimes I_Y\|_{K\otimes_{\min}   Y \to K\otimes_{\min}   Y}\to 0$ when $n\to
\infty$.

\n Hence for any $\vp>0$, for some $n$ large enough we have $\|(\sigma_n
\otimes I_Y)(v-\tilde u)\|<\vp$. Therefore,
$$\|(\sigma_n\otimes I_Y)\tilde u\| \le\|\sigma_n\|_{cb} \|v\|+\vp
<ex(Y)+\vp.$$
But on the other hand, $\tilde u-(\sigma_n\otimes I_Y)\tilde u =
((1-\sigma_n)\otimes I_Y)\tilde u \in K\otimes_{\min}  X$ since $\tilde u \in
B\otimes_{\min} X$. Hence $\hbox{dist}(\tilde u, K\otimes_{\min}   X) <ex(Y)+\vp$ and we
conclude $ex(X) <ex(Y)+\vp$.\qed

\proclaim Lemma 4. Let $X$ be an operator space. Then
$\sup\{\|T^{-1}_E\|\mid E\subset X$, $\dim E<\infty\}$ is finite iff $X$ is
exact and we have
$$ex(X) = \sup\{\|T^{-1}_E\|\mid E\subset X, \dim E<\infty\}.\leqno (8)$$

\n {\bf Proof.} Let $\lambda$ be the right side of (8). By Lemma~3 we
clearly have $\lambda \le ex(X)$ hence it suffices to show that if
$\lambda$
is finite $X$ is exact and (8) holds. Assume $\lambda$ finite. Then clearly
$T_X$ is onto. Let $q\colon \ B\otimes_{\rm min} X\to (B/K) \otimes_{\rm
min} X$ be the natural map. Consider $u$ in $\hbox{Ker}(q)$. By density
there
is a sequence $u_n$ in $B\otimes X$ such that $\|u-u_n\|<2^{-n}$. Then
$\|q(u_n)\|<2^{-n}$. By definition of $\lambda$ (since $u_n \in B\otimes E_n
$ for some finite dimensional subspace $E_n\subset X$ and 
$\|T^{-1}_{E_n}\|\le \lambda$) there is $v_n$ in $B\otimes
X$ which is a lifting of $q(u_n)$ so that $\|v_n\| <2^{-n}\lambda$ and $u_n-v_n \in K\otimes X$.
Therefore $\|u-(u_n-v_n)\|<2^{-n}+2^{-n}\lambda$ so that $u = \lim(u_n-v_n)
\in K\otimes X$. This shows that $\hbox{Ker}(q) = K\otimes_{\rm min} X$.
Thus we have showed that $\lambda<\infty$ implies
$X$ exact. By definition of $T^{-1}_X$ it is then easy to check that
$\|T^{-1}_X\|\le \lambda$.\qed

\proclaim Lemma 5. For any operator space $X$
$$ex(X)  \le \sup\{d_{SK}(E)\mid E\subset X, \dim E<\infty\}.\leqno (9)$$

\n {\bf Proof.} By the preceding lemma it suffices to show that a finite
dimensional operator space $E$ satisfies $\|T^{-1}_E\|\le d_{SK}(E)$.

\n Now consider $F\subset M_n$. By Lemma~3 and by (4)$'$ we have
$$\eqalign{\|T^{-1}_E\|  = ex(E) &\le ex(F) d_{cb}(E,F)\cr
&\le ex(M_n) d_{cb}(E,F)}$$
but trivially $ex(M_n)=1$ hence we obtain $\|T^{-1}_E\| \le d_{cb}(E,F)$
and taking the infimum over $F$, $\|T^{-1}_E\|\le d_{SK}(E)$.\qed

\n {\bf Proof of Theorem 1.} Let $E\subset X$ be finite dimensional. We will
prove
$$d_{SK}(E) \le ex(E).$$
This is the main point. To prove this claim we consider the maps $P_n\colon
\ E\to M_n$ appearing in Lemma~2. Let $E_n = P_n(E)\subset M_n$. For $n\ge
n_0$ we consider the isomorphism $u_n\colon \ E\to E_n$ obtained by
considering $P_n$ with range $E_n$ instead of $M_n$. Since $E$ is finite
dimensional $u^{-1}_n$ is $c.b.$ for each $n\ge n_0$. We claim that we have
$$\limsup_{n\to \infty} \|u^{-1}_n\|_{cb} \le ex(E).\leqno (10)$$
>From (10) it is easy to complete the proof of Theorem~1. 
Indeed, if (10) holds, we have
$$d_{SK}(E)  \le
 \limsup_{n\to \infty}\|u_n\|_{cb} \|u^{-1}_n\|_{cb} \le
ex(E).$$
By (9) we have conversely $ex(E) \le d_{SK}(E)$, whence (5). Then by (8)
$X$ is exact iff the right side of (6) is finite and (6) follows from (8).
Thus to conclude it suffices to prove our claim (10).

\n Consider $\vp_n>0$ with $\vp_n\to 0$. For each $n\ge n_0$, we choose $h_n$
in $M_{k(n)}(E)$ such that
$$\|(I_{M_{k(n)}}\otimes u_n)h_n\|_{M_{k(n)}(E_n)} = 1 \hbox{ and }
\|h_n\|_{M_{k(n)}(E)} > \|u^{-1}_n\|_{cb}-\vp_n.\leqno (11)$$
Then we form the direct sum $B_1 = \ell_\infty\{M_{k(n)}\}$ and consider
the corresponding element
$$h  =(h_n)_{n\ge n_0}\quad \hbox{in}\quad B_1 \otimes_{\rm min} E =
\ell_\infty\{M_{k(n)}(E)\}.$$
Let $K_1\subset B_1$ be the subspace formed of all the sequences $(x_n)$
with $x_n \in M_{k(n)}$ which tend to zero when $n\to \infty$. By suitably
embedding $B_1$ into $B(\ell_2)$ and $K_1$ into $K(\ell_2)$ we find that
the natural map
$$T_1\colon \ (B_1\otimes_{\rm min} E)/(K_1\otimes_{\rm min} E)\to
(B_1/K_1)\otimes_{\rm min} E$$
satisfies (note that it is an isomorphism since $\dim E<\infty$)
$$\|T^{-1}_1\| \le \|T^{-1}_E\| = ex(E).\leqno (12)$$
Let $q\colon \ B_1\otimes_{\rm min}E \to (B_1\otimes_{\rm min}
E)/(K_1\otimes_{\rm min} E)$ be the quotient mapping. Observe that we have
clearly
$$\limsup_{n\to\infty} \|h_n\| \le \|q(h)\|.\leqno (13)$$
On the other hand we have $q(h) = T^{-1}_1 T_1q(h)$ hence
$$\|q(h)\| \le \|T^{-1}_1\|\ \|T_1q(h)\|\leqno (14)$$
and since $J\colon \ E\to \ell_\infty\{E_m\}$ is a complete isometry (cf.
Lemma~2) we have
$$\|T_1q(h)\| = \|(I_{B_1/K_1}\otimes J) T_1q(h)\| _{(B_1/K_1)\otimes_{\rm
min} \ell_\infty\{E_m\}}.\leqno (15)$$
Let $q_1\colon \ B_1\otimes_{\rm min} E\to (B_1/K_1)\otimes_{\rm min}E$ be the natural map.
Clearly $q_1 = T_1q$, and the right side of (15) is the same as the norm of
the corresponding element in the space $\ell_\infty\{(B_1/K_1)\otimes_{\rm
min} E_m\}$, hence the right side of (15) is equal to
$$\eqalignno{&\sup_m\|(I_{B_1/K_1}\otimes
u_m)q_1(h)\|_{(B_1/K_1)\otimes_{\rm min} E_m}\cr
\noalign{\hbox{which is clearly}}
\le &\sup_m \limsup_{n\to \infty}\|(I_{M_{k(n)}} \otimes u_m)
(h_n)\|_{M_{k(n)}(E_m)}.}$$
For $m\le n$, we have by Lemma~2 $u_m = a_{mn}u_n$ with $\|a_{mn}\|_{cb}
\le 1$. By (11) this implies
$$\|(I_{M_{k(n)}} \otimes u_m)(h_n)\|_{M_{k(n)}(E_m)} \le \|(I_{M_{k(n)}}
\otimes u_n)(h_n)\|_{M_{k(n)}(E_n)} = 1.$$
Hence we conclude that (15) is $\le 1$. By (12), (13) and (14) we obtain
that $\limsup\limits_{n\to \infty} \|u^{-1}_n\|_{cb} = \limsup\limits_{n\to
\infty} \|h_n\| \le ex(E)$. This proves (10) and concludes the proof of
Theorem~1.\qed

 \n {\bf Remark.} The reader may have noticed that our definition of exact
 operator spaces is not the most natural extension of ``exactness'' in the
 category of operator spaces. However the more natural notion is easy to
 describe. Let us say that an operator space $X$ is
$OS$-exact if for any
 exact sequence of operator spaces (note:\ here the
morphisms are c.b. maps)
 $$\{0\}\to Y_1 \to Y_2 \to Y_3 \to \{0\}$$
 the sequence
 $$\{0\} \to Y_1 \otimes_{\rm min} X\to Y_2\otimes_{\rm min} X\to Y_3
 \otimes_{\rm min} X\to \{0\}$$
 is exact.
 
 \n Then we claim that $X$ is $OS$-exact iff there is a
constant $C$ such that
 for any finite dimensional subspace $E\subset X$, the inclusion $i_E\colon
 \ E\to X$ admits for some $n$ a factorization of the form
 $$i_E\colon \ E{\buildrel a\over \longrightarrow} M_n {\buildrel b\over
 \longrightarrow} X$$
 with $\|a\|_{cb} \|b\|_{cb} \le C$.
 
 \n Equivalently, this means that there is a net $(u_i)$ of
finite rank maps on
 $X$ of the form $u_i =b_ia_i$ with $a_i\colon \ X\to M_{n_i}$ and
 $b_i\colon \ M_{n_i}\to X$ such that
 $$\sup\|a_i\|_{cb} \|b_i\|_{cb}<\infty\quad \hbox{and}\quad
u_i(x)\to x$$
 for all $x$ in $X$.
 
 \n This result was known to E. Kirchberg and G. Vaillant.
It can be proved as follows. First if
 $X$ is $OS$-exact, it is a fortiori exact in the above sense so that by
 Theorem~1 there is a constant $C_1$ such that $d_{SK}(E)\le C_1$ for all
 finite dimensional subspaces $E\subset X$.
 
 \n Secondly, if $X$ is $OS$-exact, there is clearly a
constant $C_2$ such that
 for any pair of finite dimensional operator spaces $E_1\subset E_2$ (so
 that $E^\bot_1 \subset E^*_2$) we have an isomorphism
 $$T\colon \ (E^*_2\otimes_{\rm min} X)/(E^\bot_1
\otimes_{\rm min}
 X)\longrightarrow (E^*_2/E^\bot_1) \otimes_{\rm min} X$$
 such that $\|T^{-1}\|\le C_2$.
 
 \n In other words, since $E^*_2\otimes_{\rm min} X =
cb(E_2,X)$ we have an
 extension property associated to the following diagram:
 $$\matrix{E_2\cr
 &\quad \searrow \widetilde v\cr
 \cup\cr
 E_1&{\buildrel v\over \longrightarrow}&X\cr}$$
 More precisely, for any $v\colon \ E_1\to X$ there is an extension
 $\widetilde v\colon\ E_2\to X$ such that $\|\widetilde v\|_{cb}\le
 C_2\|v\|_{cb}$. Consider now an arbitrary finite dimensional subspace
 $E\subset X$. Let $\vp>0$. Consider $E_1\subset M_n$ such that there is an
 isomorphism $u\colon \ E_1\to E$ with $\|u\|_{cb} \|u^{-1}\|_{cb}\le
 d_{SK}(E) +\vp \le C_1+\vp$.
 
 \n Using the preceding extension property (with $E_2=M_n$)
we find an operator
  $b\colon \ M_n\to X$ extending $u$ and such that
$\|b\|_{cb}\le
 C_2\|u\|_{cb}$. Let $a\colon \ E\to M_n$ be the operator $u^{-1}$
 considered as acting into $M_n$. Then $i_E=ba$ and
 $$\|a\|_{cb}\|b\|_{cb}
 \le C_2\|u\|_{cb} \|u^{-1}\|_{cb} \le C_2(C_1+\vp).$$
 This proves our claim. Note in particular
 that if $X$ is a $C^*$-algebra, it is $OS$-exact
 iff it is nuclear, by [P1, Remark before Theorem~2.10].

 In the category of Banach spaces one can define a similar notion of
 exactness using the {\it injective\/} tensor product instead of the minimal
 one. Then a Banach space is ``exact'' iff it is a ${\cal L}_\infty$-space
 in the sense of [LR]. We refer the reader to [LR, Theorem~4.1 and subsequent
 Remark]. 
\vskip24pt
\vfill\eject
\n {\bf \S 2. Ultraproducts.}

The notion of exactness for operator spaces is closely connected to a commutation property
involving ultraproducts. To explain this let us recall a few facts about
ultraproducts. Let $({ F}_i)_{i\in I}$ be a family of operator spaces
and let
${\cal U}$ be a nontrivial ultrafilter on $I$. 
We denote by $\widehat { 
F} = \Pi { F}_i/{\cal U}$ the associated ultraproduct in the category of
Banach spaces (cf. e.g. [Hei]). Recall that if   $dim(F_i)=n$ for all $i$ in $I$,
 then
the ultraproduct $\widehat { 
F}$ clearly also is $n$-dimensional.

Clearly $\widehat { F}$ can be equipped with an operator space structure
by defining
$$M_n(\widehat{ F}) = \Pi M_n({ F}_i)/{\cal U}.\leqno (16)$$
It is easy to check that Ruan's axioms  [Ru] are satisfied so that
$\widehat { F}$ with the matricial structure (16) is an operator space.
Alternatively, one may view ${ F}_i$ as embedded into $B(H_i)$ ($H_i$
Hilbert) and observe that $\widehat{ F} \subset \Pi B(H_i)/{\cal U}$.
Since $C^*$-algebras are stable by ultraproduct
 we obtain $\widehat { 
F}$ embedded in a $C^*$-algebra. It is easy to see that the resulting
operator space structure is the same as the one defined by (16). Note that
(16) can be written as a commutation property between ultraproducts and
 the minimal tensor product, as follows
$$M_n \otimes_{\rm min}[\Pi { F}_i/{\cal U}] = \Pi[M_n\otimes_{\rm
min} { F}_i]/{\cal U}.\leqno (17)$$
It is natural to wonder which operator spaces $E$ can be substituted to
$M_n$
in this identity (17). It turns out that this property is closely related
to the invariant $d_{SK}(E)$, as we will now show.

We first observe that there is for any finite dimensional operator space
$E$ a canonical map
$$v_E\colon \ \Pi(E\otimes_{\rm min}E_i)/{\cal U} \to E\otimes_{\rm min}
\widehat E\leqno (18)$$
with $\|v_E\|\le 1$. Indeed, we clearly have a norm one mapping
$$E\otimes_{\rm min} \ell_\infty\{E_i\} \to E\otimes_{\min}  \widehat E\leqno(18)'$$
but if $E$ is finite dimensional $E\otimes_{\rm min} \ell_\infty\{E_i\} =
\ell_\infty\{E\otimes_{\rm min}E_i\}$ and the map (18)' vanishes on the
subspace of elements with ${\cal U}$ limit zero. Hence, after passing to the quotient by the kernel of (18)',
 we find the map (18)
with norm $\le 1$. More generally (recall the isometric identity
$F^*\otimes_{\rm min} E = cb(F,E)$, cf. [BP, ER1]) if $(E_i)_{i\in I}$
(resp. $(F_i)_{i\in I}$) is a family of $n$-dimensional (resp.
$m$-dimensional operator spaces), we clearly have a norm one canonical map
$$\Pi cb(E_i,F_i)/{\cal U} \to cb(\widehat E, \widehat F),\leqno (18)''$$
where $\widehat E  = \Pi E_i/{\cal U}$ and $\widehat F = \Pi F_i/{\cal U}$.

\proclaim Proposition 6. Let $E$ be a finite dimensional operator space and
let $C\ge 1$ be a constant. The following are equivalent.\medskip
\item{(i)} $d_{SK}(E) \le C$.
\item{(ii)} For all ultraproducts $\widehat F = \Pi F_i/{\cal U}$ the
canonical isomorphism (which has norm $\le 1$)
$$v_E\colon \ \Pi(E\otimes_{\rm min} F_i)/{\cal U}\to E\otimes_{\rm min}
(\Pi F_i/{\cal U})$$
satisfies $\|v^{-1}_E\|\le C$.
\item{(iii)} Same as (ii) but with all ultraproducts $(F_i)_{i\in I}$ on a
countable set and such that $\sup\limits_{i\in I} \dim F_i \le \dim
E$.\medskip

\pf First observe that if $G\subset F$ are operator spaces then we have
isometric embeddings
$$G\otimes_{\rm min} \widehat F\to F\otimes_{\rm min} \widehat F$$
and
$$\Pi(G\otimes_{\rm min} F_i)/{\cal U}\to \Pi(F\otimes_{\rm min}
F_i)/{\cal U}.$$
Therefore in the finite dimensional case we have clearly $\|v^{-1}_G\| \le
\|v^{-1}_F\|$.

\n (i) $\Rightarrow$ (ii): Assume (i). Then consider $G\subset M_n$
isomorphic to $E$. We have clearly $\|v^{-1}_{M_n}\| =1$ by (16), hence we
can write
$$\eqalign{\|v^{-1}_E\| \le d_{cb}(E,G)\|v^{-1}_G\| &\le
d_{cb}(E,G)\|v^{-1}_{M_n}\|\cr
&\le d_{cb}(E,G)}$$
hence $\|v^{-1}_E\| \le d_{SK}(E)$, whence (ii).

\n  (ii) $\Rightarrow$ (iii) is trivial.

\n (iii) $\Rightarrow$ (i): This is proved by an argument similar to the
proof of (10) in Theorem~1. We merely outline the argument. Let $u_n\colon
\ E\to E_n = P_n(E)$ be given by Lemma~2, as in the above proof of (10). For $n\ge n_0$ we consider
$u^{-1}_n\colon \ E_n\to E$ and we identify $u^{-1}_n$ with an element of
$E^*_n \otimes E$. Recall $\|u^{-1}_n\|_{E^*_n\otimes _{\rm min} E} =
\|u^{-1}_n\|_{cb(E_n,E)}$. Then
$$\|(u^{-1}_n)_n\|_{\Pi(E^*_n\otimes_{\rm min}E)/{\cal U}} = \lim_{\cal
U}\|u^{-1}_n\|_{cb}$$
and on the other hand since $J$ is a complete isometry and since we have the
monotonicity property (7) we have
$$\eqalign{\|(u^{-1}_n)_n\|_{(\Pi E^*_n/{\cal
U})\otimes_{\rm min} E} &=\sup_m \|(P_m u^{-1}_n)_n\|_{(\Pi
E^*_n/{\cal U})\otimes_{\rm min}M_m}\cr  &=\sup_m
\lim_{n,{\cal U}}\|P_m u^{-1}_n\|_{ E^*_n\otimes_{\rm min}
M_m} \cr &\le \sup_m \lim_{n,{\cal U}}\| a_{mn}\|_{cb}
 \le 1.}$$
Hence (iii) implies $\lim\limits_{\cal U}\|u^{-1}_n\|_{cb}\le C$ if ${\cal
U}$ is any nontrivial ultrafilter on $\bf N$, and  we conclude  $d_{SK}(E)
\le \lim\limits_{\cal U} \|u_n\|_{cb} \|u^{-1}_n\|_{cb} \le C$.\qed
\vskip24pt

\n {\bf \S 3. How large can $d_{SK}(E)$ be?}

We now wish to produce finite dimensional operator spaces $E$ with
$d_{SK}(E)$ as large as possible. It follows from 
Theorem 9.6 in [P1] that for any $n$-dimensional operator
space $E$ we have
$$d_{SK}(E) \le \sqrt n.$$
We will show that this upper bound in general cannot be improved,
 at least asymptotically, when $n$ goes to infinity.
We will consider the space $\ell^n_\infty$ with its natural operator space
structure. We will denote by $E^n_1$ the dual in the category of Banach
spaces, so that as a Banach space $E^n_1$ is the usual space $\ell^n_1$,
however it is embedded into $B(H)$ in such a way that the canonical basis
$e_1,\ldots, e_n$ of $E^n_1$ satisfies for all $a_1,\ldots, a_n$ in
$B=B(\ell_2)$.
$$\left\|\sum^n_1 e_i\otimes a_i\right\|_{E^n_1\otimes_{\rm min}B} =
\sup_{\scriptstyle u_i\atop \scriptstyle {\rm unitary}} \left\|\sum^n_1
u_i\otimes a_i\right\|_{B\otimes_{\rm min}B} \leqno (19)$$
where the supremum runs over all unitary operators $u_i$ in $B$.

Another remarkable representation of $E^n_1$ appears if we consider the
full $C^*$-algebra $C^*({\bf F}_n)$ of the free group with $n$ generators.
If we denote by $\delta_1,\ldots, \delta_n$ the generators of ${\bf F}_n$
viewed as unitary operators in $C^*({\bf F}_n)$ in the usual way then
$\|\sum \delta_i\otimes a_i\|_{C^*({\bf F}_n)\otimes_{\rm min}
B} $ is equal to (19), which shows that the map $u\colon \ E^n_1\to
\hbox{span}(\delta_i)$ which takes $e_i$ to $\delta_i$ is a complete
isometry. Our main result is the following.

\proclaim Theorem 7. For all $n\ge 2$
$$d_{SK}(E^n_1) \ge {n\over 2\sqrt{n-1}}.$$
Hence in particular $d_{SK}(E^n_1)>1$ for all $n\ge 3$.

\n {\bf Remark.} It is easy to check (this was pointed out to me by
Paulsen) that $d_{SK}(E^n_1)=1$ for $n=2$. Indeed, in that case (19)
becomes (after multiplication by $u^{-1}_1$)
$$\|e_1\otimes a_1+e_2\otimes a_2\| = \sup_{u \ {\rm unitary}} \|I\otimes
a_1 + u\otimes a_2\|$$
and since $(I,u)$ generate a commutative $C^*$-algebra, this is the same as
$$\sup_{\scriptstyle |z|=1\atop\scriptstyle z\in {\bf C}} \|a_1+za_2\|$$
which shows that $E^2_1$ is completely isometric to the span  of $\{1,
e^{it}\}$ in $C({\bf T})$. Therefore (since $C({\bf T})$ is nuclear)
$d_{SK}(E^2_1)=1$. We will use an idea similar to Wassermann's argument in
[W1]:\ we consider the direct sum $M = M_1 \oplus M_2\oplus\cdots$, or
equivalently $M=\ell_\infty\{M_n\}$ in our previous notation, and we denote
by $I_{\cal U}$ the set
$$I_{\cal U} = \{(x_\alpha)\in M\mid \lim_{\cal U}
\tau_\alpha(x^*_\alpha x_\alpha)=0\}$$
where $\tau_\alpha$ is the normalized trace on $M_\alpha$ and where ${\cal
U}$ is a nontrivial ultrafilter on $\bf N$. Then the group von~Neumann
algebra $VN({\bf F}_n)$ is isomorphic to a von~Neumann subalgebra of the
quotient $N=M/I_{\cal U}$. It is well known that $N$ is a finite
von~Neumann
algebra with normalized trace $\tau$ given by $\tau(x) = \lim\limits_{\cal
U} \tau_\alpha(x_\alpha)$ where $x$ denotes the equivalence class  in $N$
of $(x_\alpha)$. Let us denote by
$$\Phi\colon \ M\to M/I_{\cal U}$$
the quotient mapping.

In the sequel, we will make use of the operator space version of the
projective tensor product introduced in [ER5]. However, to facilitate the
task of the reader, we include in the next few lines the simple facts that
we use with indication of proof. Consider a finite dimensional algebra
$M_k$ equipped with the normalized trace which we denote by $\tau$.
We denote by $L_1(\tau)$ the space $M_k$ equipped with the
norm $\|x\|_1 =\tau(|x|)$. Let $E$ be an operator
space. Since $M_k = L_1(\tau)^*$ we have $(L_1(\tau)\otimes E)^* =
M_k(E^*)$. We then denote by $L_1(\tau)\otimes_\wedge E$ the space
$L_1(\tau)\otimes E$ equipped with the norm induced on $L_1(\tau)\otimes E$
by $M_k(E^*)^*$. We will use the following two facts which are easy to
check:

(a)\ If $F\subset E$ (completely isometric embedding) then $L_1(\tau)
\otimes_\wedge F\subset L_1(\tau) \otimes_\wedge E$ (completely isometric
embedding).

(b) \ If $E=E^n_1$ and $e_1,e_2,\ldots, e_n$ is the canonical basis of
$E^1_n = {\ell^{n}_\infty}^*$, then for any $x_1,\ldots, x_n$ in $L_1(\tau)$
we have
$$\left\|\sum^n_1 x_i\otimes e_i\right\|_{L_1(\tau)\otimes_\wedge E^1_n} =
\sum^n_1\|x_i\|_{L_1(\tau)}.$$

(c)\ We have a norm one inclusion
$$M_k(E)\to L_1(\tau) \otimes_\wedge E.$$
These facts can be checked as follows:\medskip

\item{(a)} follows by duality from the isometric identity
$$M_k(F^*) = M_k(E^*/F^\bot)=M_k(E^*)/M_k(F^\bot).$$
\item{(b)} follows again by duality from the identity
$$M_k(\ell^n_\infty) = \ell^n_\infty(M_k).$$
\item{(c)} follows from the inequality $\forall\ \xi_{ij}\in E^*$ $\forall
\ x_{ij}\in E$

$$k^{-1}\left|\sum_{ij\le k}\xi_{ij}(x_{ij})\right| \le
\|(x_{ij})\|_{M_k(E)} \|(\xi_{ij})\|_{M_k(E^*)}.\leqno(d)$$
The latter inequality can be checked using the factorization of $c.b.$
maps (cf. [Pa1, p. 100]) since 
$\|(\xi_{ij})\|_{M_n(E^*)}=\|(\xi_{ij})\|_{cb(E,M_n)}$:
\ if $\|(\xi_{ij})\|_{M_n(E^*)}\le 1$ then we can write
$\xi_{ij}(x) =
\langle\pi(x)x_j,y_i\rangle$ with $\pi\colon \ E\to B(H)$ restriction of a
representation and with $x_j,y_i \in H$ such that $\|\sum \alpha_jx_j\|\le
1$ and $\|\sum \alpha_iy_i\|\le 1$ whenever $\sum|\alpha_j|^2 \le 1$. This
implies $(\sum\|x_j\|^2 \sum\|y_i\|^2)^{1/2} \le k$ whence
$$\left|\sum \xi_{ij}(x_{ij})\right|  = \left|\sum \langle
\pi(x_{ij})x_j,y_i\rangle\right| \le k \|(x_{ij})\|_{M_k(E)}$$
which proves the inequality $(d)$.

 We denote below by $ C^*_\lambda({\bf F}_n)$
the reduced $C^*$-algebra associated to the left regular representation
for  the free group ${\bf F}_n$ with $n$ generators. Then the key result for our subsequent 
estimates can be stated as follows.

\proclaim Theorem 8. Fix $n\ge 2$. There is a family of unitary matrices
$(u^\alpha_i)$ with $u^\alpha_i \in M_\alpha$, $(i=1,\ldots, n$, $\alpha\in
{\bf N})$ and a nontrivial ultrafilter ${\cal U}$ on $\bf N$ such that for
all $m\ge 1$ and all $x_1,\ldots, x_n$ in $M_m$ we have
$$\lim_{\cal U} \left\|\sum^n_{i=1} u^\alpha_i \otimes
x_i\right\|_{L_1(\tau_\alpha)\otimes_\wedge M_m} \le\left\|\sum^n_1
\lambda(g_i)\otimes x_i\right\|_{ C^*_\lambda({\bf F}_n) \otimes_{\rm min}
M_m}\leqno (20)$$

\n {\bf Remark.} By results included in [HP], the right side of (20) is
$$\le 2 \max\left\{\left\|\sum\limits^n_1
x^*_ix_i\right\|^{1/2}\right. ,
\left.\left\|\sum\limits^n_1
x_ix^*_i\right\|^{1/2}\right\}.$$

\n {\bf Proof of Theorem  8.} As explained in [W1], for each $i$ there is a
unitary $u_i = (u^\alpha_i)_{\alpha \in {\bf N}}$ in $M$ such that
$\Phi(u_i) = \lambda(g_i)$. Let us denote
$$b = \sum \lambda(g_i) \otimes x_i\in C^*_\lambda({\bf F}_n)\otimes_{\rm min}
M_m.$$
Then $\sum u_i\otimes x_i$ is a lifting of $b$ in $M\otimes M_m$ and
$M_m(VN({\bf F}_n))$ embeds isometrically (see [W1]) into $M_m(M/I_{\cal
U})$ or equivalently into $M_m(M)/M_m(I_{\cal U})$. It follows that we have
$$\|b\| = \inf\left\{\left\| \sum^m_{i=1} x_i\otimes u_i
+\gamma\right\|_{M_m(M)} \mid \gamma \in M_m(I_{\cal U})\right\}.$$
Hence there is a sequence $(\gamma^\alpha)_{\gamma\in {\bf N}}$ with
$\gamma^\alpha \in M_m(M_\alpha)$ satisfying
$$\forall\ i,j\le m\qquad \lim_{\scriptstyle \alpha\to \infty \atop
\scriptstyle {\cal U}} \tau_\alpha(\gamma^{\alpha
*}_{ij}\gamma^\alpha_{ij}) = 0$$
and such that
$$\lim_{\cal U} \left\|\sum^n_{i=1} x_i\otimes u^\alpha_i
+\gamma^\alpha\right\|_{M_m(M_\alpha)} \le \|b\|.\leqno (21)$$
Now observe that $M_m(M_\alpha) = M_\alpha(M_m)$. We will use the norm one
inclusion (see fact (c) above) $M_\alpha(M_m)\longmapsto L_1(\tau_\alpha)
\otimes_\wedge M_m$. Note that the inclusion $L_2(\tau_\alpha)\to
L_1(\tau_\alpha)$ has norm $\le 1$ so that
$$\lim_{\scriptstyle \alpha\to \infty \atop
\scriptstyle {\cal U}} \|\gamma^\alpha\|_{L^1(\tau_\alpha)\otimes_\wedge
M_m} =0.$$
Therefore (21) yields
$$\lim_{\cal U} \left\|\sum\nolimits^n_1 u^\alpha_i \otimes
x_i\right\|_{L^1(\tau_\alpha)\otimes_\wedge M_m} \le \|b\|,\leqno (22)$$
which is the announced inequality.\qed

To prove Theorem~7, we will use the following lemma. 

\proclaim Lemma 9. Consider the operator
$$T_n\colon \ \ell^n_\infty\to C^*_\lambda({\bf F}_n)$$
defined by $T_n(\alpha_1,\ldots, \alpha_n) = \sum\limits^n_1
\alpha_i\lambda(g_i)$. Then for any $m$, any subspace $F\subset M_m$ and
any factorization
$$\matrix{&F^*\cr
&a\nearrow \qquad \qquad b\searrow\cr
\ell^n_\infty&{\buildrel T_n \over
{\hbox to 35pt{\rightarrowfill}}}&C^*_\lambda({\bf F}_n)\cr}$$
with $T_n = ba$ we have
$$n \le \|a\|_{cb} \|b\|_{cb}.\leqno (23)$$

\n {\bf Proof.} Consider $a,b$ as above. We identify $b\colon\ F^*\to
C^*_\lambda({\bf F}_n)$ with an element of \break $F\otimes_{\rm min}
C^*_\lambda({\bf F}_n)$. Then we can write $b = \sum\limits^n_{i=1}
x_i\otimes \lambda(g_i)$ with $x_i\in F$ such that $a^*(x_i) =e_i$. (Recall
that $e_i$ is the canonical basis of $E^n_1 = (\ell^n_\infty)^*$.)

\n Now by Theorem~8 we have
$$\lim_{\cal U} \left\|\sum\nolimits^n_1 u^\alpha_i \otimes
x_i\right\|_{L_1(\tau_\alpha)\otimes_\wedge M_m} \le \|b\|_{cb}.$$
By fact (a) recalled above, this implies
$$\lim_{\cal U} \left\|\sum\nolimits^n_1 u^\alpha_i\otimes
x_i\right\|_{L_1(\tau_\alpha)\otimes_\wedge F} \le \|b\|_{cb},$$
 hence since $\|a^*\|_{cb} =
\|a\|_{cb}$
$$\lim_{\cal U} \left\|\sum\nolimits^n_1 u^\alpha_i \otimes
a^*(x_i)\right\|_{L_1(\tau_\alpha)\otimes_\wedge E^n_1} \le \|a\|_{cb}
\|b\|_{cb}.$$
This gives the conclusion since $a^*(x_i) = e_i$ and by fact (b)
$$\left\|\sum\nolimits^n_1 u^\alpha_i \otimes
e_i\right\|_{L_1(\tau_\alpha)\otimes_\wedge E^n_1} =
\sum\nolimits^n_1\|u^\alpha_i\|_{L_1(\tau_\alpha)} = n.$$
\qed

\n {\bf Remark.} The same operator
$T_n$ as in Lemma 9  was already 
considered  in [H]. By [H, Lemma 2.5] we have
$\|T_n\|_{dec}=n$, but this does not seem related to (23).   

\n {\bf Proof of Theorem 7.} By [AO] (see [H] for more details) we have
$\|T_n\|_{cb} \le 2\sqrt{n-1}$, hence for any $F\subset M_m$ we can write
by (23)
$$\leqalignno{n &\le \|T_n\|_{cb} d_{cb}(F^*, \ell^n_\infty)&(24)\cr
n &\le \|T_n\|_{cb} d_{cb}(F^*, E_n^{\lambda})&(25)}$$
where $E_n^{\lambda} = \hbox{span}\{\lambda(g_i)\mid i=1,\ldots, n\}$ in
$C^*_\lambda({\bf F}_n)$. Since $d_{cb}(F^*, \ell^n_\infty) =
d_{cb}(F,{\ell^{n}_\infty}^*)$ and ${\ell^{n}_\infty}^* =
E^n_1$ we obtain $$n(2\sqrt{n-1})^{-1} \le
d_{cb}(F,E^n_1)$$ so that Theorem~7 follows.

\def\T{\theta}
\n {\bf Remark.} By the same argument we have
$$d_{SK}({E_n^{\lambda}}^*) \ge n(2\sqrt{n-1})^{-1}.$$
Here again this is $>1$ if $n\ge 3$ but $d_{SK}({E_n^{\lambda}}^*)=1$ if $n=2$ for the
same reason as above for $E^n_1$. We can also derive an estimate for the
$n$ dimensional operator Hilbert space which is denoted by $OH_n$. This
space was introduced in [P1] to which we refer for more details. It
is isometric to $\ell^n_2$ and has an orthonormal basis $(\T_i)_{i\le n}$
such that for all $a_1,\ldots, a_n$ in $B$ we have
$$\left\|\sum^n_1 \T_i\otimes a_i\right\|_{OH_n\otimes_{\rm min} B} =
\left\|\sum^n_1 a_i \otimes \bar a_i\right\|^{1/2}_{B\otimes_{\rm
min}\overline B}.\leqno (26)$$

\proclaim Corollary 10. For each $n\ge 2$, we have $$d_{SK}(OH_n) \ge
[n(2\sqrt{n-1})^{-1}]^{1/2}.\leqno(27)$$

\n {\bf Proof.} By (26) we have (we denote simply by $\| \ \|$ the minimal
tensor norm everywhere)
$$\left\|\sum \T_i\otimes \lambda(g_i)\right\| = \left\|\sum \lambda(g_i)
\otimes \overline{\lambda(g_i)}\right\|^{1/2}$$
and by [AO]  we have
$$\left\|\sum^n_1 \lambda(g_i) \otimes \overline{\lambda(g_i)}\right\| =
2\sqrt{n-1}.$$
Hence we have a factorization of $T_n$ of the form $\ell^n_\infty
{\buildrel a\over \longrightarrow} OH_n {\buildrel b \over \longrightarrow}
E_n$ with $\|b\|^2_{cb} \le 2\sqrt{n-1}$. 
On the other hand $$\|a\|_{cb} =
\|a^*\|_{cb} =
\left\|\sum\limits^n_1 e_i\otimes \bar
e_i\right\|_{E^n_1\otimes_{\min} \overline{E^n_1}}^{1/2}
\le n^{1/2}$$ hence we have by (23) $$n \le d_{SK}(OH_n)
\|a\|_{cb} \|b\|_{cb}$$
which implies (27).\qed

\n {\bf Remark.} It is easy to verify that these estimates are asymptotically
best possible. More precisely, we have with the notation
of [P1],  $d_{SK}(E^n_1) \le d_{cb}(E^n_1,R_n)\le
{\sqrt{n}}.$ Similarly, $d_{SK}({E_n^{\lambda}}^*) \le
d_{cb}({E_n^{\lambda}}^*,R_n)\le {\sqrt{n}},$ and finally
 $d_{SK}(OH_n) \le d_{cb}(OH_n,R_n\cap C_n ) \le
{{n}^{1/4}}.$

\n {\bf Remark.} It is natural to raise the following question:\ Is there a
function $n\to f(n)$ and a constant $C$ such that for any
$E$ in $OS_n$ satisfying $d_{SK}(E)=1$, there is a
subspace $F\subset M_m$ with $m\le f(n)$ and $d_{cb}(F,E)
\le C$?  It is easy to derive from the preceding construction
that the answer is negative (contrary to the commutative
case with $\ell^n_\infty$ in the place of $M_n$). A
negative answer can also be derived easily from the fact
(due to Szankowski [Sz]) that the space of compact
operators on $\ell_2$ fails the {\it uniform\/}
approximation property.

\n {\bf Remark.} It is rather natural to introduce the following quantity
for $E$ in $OS_n$.
$$d_{QSK}(E) = \inf\{d_{cb}(E,F)\}$$
where the infimum runs over all the spaces $F$ which are quotient of a
subspace of $K$, i.e. there are $S_2\subset S_1\subset K$ such that $F =
S_1/ S_2$. Clearly $d_{QSK}(E) \le d_{SK}(E)$. We do not
know much about this new parameter. Note however that in
view of the lifting property of $E^n_1$ we have by
Theorem~7 $$d_{QSK}(E^n_1) \ge {n\over 2\sqrt{n-1}}.$$
Moreover, it can be shown  that
$OH_n$ embeds completely isometrically into the direct sum
$L_\infty (R_n)\oplus L_\infty (C_n)$, and a fortiori
into $L_\infty (M_{2n})$, so
that $d_{QSK}(OH_n) =1$, indeed this follows from the
identity $OH_n=(R_n,C_n)_{1/2}$ proved in  [P1, Corollary
2.6]. Note that for the quotient of a subspace case, the
identity between (1) and (1)' might no longer hold, so
that it might be necessary to distinguish between the quotients of a subspace
of $K$ and those of $M_n$ for some $n\ge 1$.

\vskip24pt
\n {\bf \S 4. On the metric space of all $n$-dimensional operator spaces.}

It will be convenient in the sequel to record in the next statement
  several elementary facts
on ultraproducts.

\proclaim Proposition 11. Let $E$ and $F$
be two $n$-dimensional operator spaces.
Let $(E_i)_{i\in I}$ and $(F_i)_{i\in I}$
be two
families of $n$
dimensional operator spaces. Let ${\cal U}$ be an ultrafilter
on $I$ and let $\widehat E$ ,\ $\widehat F$ 
be the corresponding
ultraproducts.  
\item {(i)}\ We have
$$d_{cb}(\widehat E, \widehat F)  \le \lim_{\cal U} d_{cb}(E_i,F_i).\leqno (28)$$
\n Moreover, if $d_{cb}(E_i,F_i)\to 1$ then $\widehat E$ is
completely isometric to $\widehat F$.
\item {(ii)}\ If $F_i =F$ for all $i\in I$ then $\widehat F$ is completely
isometric to $F$.
\item {(iii)}\ If $d_{cb}(E_i,F)\to 1$ then $\widehat E$ is completely
isometric to $F$.
\item {(iv)}\ If $d_{cb}(E,F)= 1$, then $E$ and $F$ are completely
isometric.

\pf \ We have isomorphisms $u_i\colon \ E_i\to F_i$ such that $v_i =
u^{-1}_i$ satisfies $\|u_i\|_{cb}\le 1$ and $\lim\limits_{\cal U}
\|v_i\|_{cb} = \lim\limits_{\cal U} d_{cb}(E_i,F_i)$. Then 
the maps $\hat
u\colon \ \widehat E\to \widehat F$ and $\hat v\colon \ \widehat F\to
\widehat E$  (associated respectively to $(u_i)$ and $(v_i)$) are inverse of each other and satisfy $\|\hat u\|_{cb}\|\hat
v\|_{cb} \le \lim\limits_{\cal U} d_{cb}(E_i,F_i)$.
 Hence we have (28).
The preceding also shows that $\hat u$ is a complete isometry between
$\widehat E$ and $\widehat F$ when $\lim\limits_{\cal U} d_{cb}(E_i,F_i) =
1$, whence (i). But on the other hand it is easy to check that, if $F_i=F$ for all $i$, then $\widehat F$ and $F$ are
completely isometric via the map $T\colon \ \widehat F\to F$ defined by
$T(x) = \lim\limits_{\cal U}x_i$ if $x$ is the equivalence class of $(x_i)$
modulo ${\cal U}$. This justifies (ii). Then (iii) is clear.
Finally, taking $E_i=E$ and $F_i=F$ for all $i$ in what precedes,
we obtain (iv). (Actually, (iv) is clear, by a direct argument
based on the compactness of the unit ball of $cb(E,F)$.)\qed

In this section, we will include some remarks on the set
$OS_n$ formed of all $n$ dimensional operator spaces. More precisely
$OS_n$
is the set of all equivalence classes when we identify two spaces if they
are completely isometric. We equip $OS_n$ with the metric
$$\delta(E,F) = \hbox{Log } d_{cb}(E,F).$$
This is the analogue for operator spaces of the classical ``Banach-Mazur
compactum'' formed of all $n$ dimensional normed spaces equipped with the
Banach-Mazur distance.

However, contrary to the Banach space situation the metric space $OS_n$ is
not compact. The next result was known, it was mentioned to me by Kirchberg.

\proclaim Proposition 12. The set $OS_n$ equipped with the metric
$\hbox{Log } d_{cb}$ is a complete metric space, but it is not compact at
least if $n\ge 3$.

\pf We sketch a proof using ultrafilters. Let $(E_i)$ be a Cauchy sequence
in $OS_n$. Let $\widehat E$ be an ultraproduct associated to a nontrivial
ultrafilter ${\cal U}$ on $\bf N$. Then by the Cauchy condition for each
$\vp>0$ there is an $i_0$ such that for all $i,j>i_0$ we have
$$\hbox{Log } d_{cb}(E_i,E_j) \le \vp.$$
By Proposition 11 this implies $\forall\ i>i_0$
$$\hbox{Log } d_{cb}(\widehat E,E_i)\le \vp$$
hence $E_i\to \widehat E$ when $i\to \infty$ and $OS_n$ is complete.

We now show that $OS_n$ is not compact by exhibiting a sequence without
converging subsequences if $n\ge 3$. Consider any space $E_0$ in $OS_n$
such that
$$d_{SK}(E_0)>1.\leqno (29)$$
We know that such spaces exist if $n\ge 3$ by Theorem~7 and Corollary~10.

\n By Lemma~2, we can find a sequence of spaces $E_i\subset M_i$ with $\dim
E_i = n$ such that for any nontrivial ultrafilter on $\bf N$ the
ultraproduct $\widehat E = \Pi E_i/{\cal U}$ is completely isometric to
$E_0$. (Indeed, this is clear by (7).) Assume that some subsequence of
$E_i$ converges in $OS_n$. Then, by Proposition~11 (ii), its limit must be
$\widehat E$ which is the same as $E_0$. In other words the subsequence can
only converge to $E_0$ but on the other hand by (29) we have (since
$E_i\subset M_i$)
$$\forall\ i\in I\qquad 1< d_{SK}(E_0) \le d_{cb}(E_0,E_i),$$
which is the desired contradiction.\qed

\n {\bf Remark.} It is well known that every separable Banach space $E$
embeds isometrically into the space of all continuous functions on the unit
ball of $X^*$ equipped with the weak$^*$-topology.
 Let us denote simply by
${\cal C}$ the latter space, and let $k\colon \ E\to {\cal C}$ be the isometric
embedding. Given $P_n$ as in Lemma~2, we can introduce $\widetilde
P_n\colon \ E\to {\cal C}\oplus_\infty M_n$ by setting $\widetilde P_n(x) =
(k(x), P_n(x))$. Then each $\widetilde P_n$ is an
 {\it isometric\/}
isomorphism of $E$ into $\widetilde E_n \subset {\cal C}\break \oplus_\infty M_n$.
Moreover the embedding $\widetilde J\colon \ E\to \ell_\infty\{\widetilde
E_n\}$ is a completely isometric embedding with the same properties as in
Lemma~2. Finally the $C^*$-algebras ${\cal C}\oplus_\infty M_n$ are nuclear. Using
this it is easy to modify the preceding reasoning, replacing $E_n$ by
$\widetilde E_n$ (note $d_{SK}(\widetilde E_n) =1$ since ${\cal C}\oplus_\infty
M_n$ is nuclear) and to demonstrate the following

\proclaim Corollary 13. Let $E_0$ be any $n$ dimensional operator space
such that $d_{SK}(E_0)>1$. Then the (closed) subset of $OS_n$ formed of all the
spaces isometric to $E_0$ is not compact.

\n {\bf Remark.} Actually we obtain a sequence $E_i$ of spaces each isometric to $E_0$ and
such that $E_0$ is completely isometric to $\widehat E =\Pi E_i/{\cal U}$
but $d_{cb}(E_i,E_0)\not\to 0$. More precisely (in answer to a question of S.~Szarek), the
preceding argument shows that the ``metric entropy'' of $OS_n$ is quite
large in the following sense:\ Let $\delta = d_{SK}(E_0)$ with $E_0$ as in
Corollary~13. Then for any $\vp>0$ there is a sequence $E_i$ in $OS_n$ such
that
$$d_{cb}(E_i,E_j) > \delta-\vp\quad \hbox{for any}\quad i\ne j.$$
(Moreover, $E_i$ is isometric to $E_0$ and the ultraproduct $\widehat E =
\Pi E_i/{\cal U}$ is completely isometric to $E_0$). This suggests the
following question:\ does there exist such a sequence if $\delta$ is equal
to the diameter of the set $OS_n$ (or of the subset of $OS_n$ formed of all
the spaces which are isometric to $E_0$)? If not, what is
the critical value of $\delta$?

Corollary 13 is of course particularly striking in the case $E_0 = OH_n$,
it shows that the set of all possible operator space structures on the
Euclidean space $\ell^n_2$ is very large. We refer to
 [Pa2] for more information of the latter set.

These results lead to the following question (due to Kirchberg).\medskip

\n {\bf Problem.} Is the metric space $OS_n$ separable?

\n Equivalently, is there a separable operator space $X$
such that for any
$\epsilon >0$ and and any $n$-dimensional operator space $E$,
there is a subspace $F\subset X$ such that $d_{cb}(E,F)<1+\epsilon$?

\n More generally, let $E$ be an $n$-dimensional normed
space. Let $OS_n(E)$ be the subset of $OS_n$ formed of all
the spaces which are isometric to $E$. For which spaces $E$
is the space   $OS_n(E)$ compact or separable? 
Note 
that $OS_n(E)$ can be a singleton,
 this happens in the $2$-dimensional
case if
$E=\ell_1^2$ or $E=\ell_\infty^2$, however it never happens
if $dim(E)\ge 5$ (see [Pa2, Theorem 2.13]).

We will now characterize the spaces $E_0$ for which the conclusion of
Corollary~13 holds. We will use the following simple observation.

\proclaim Lemma 14. Fix $n\ge 1$. Let $\widehat E = \Pi E_i/{\cal U}$ be
an ultraproduct of $n$-dimensional spaces. Then $(\widehat E)^* = \Pi
E^*_i/{\cal U}$ completely isometrically.

\pf Let $F_i$ be a family of $m$-dimensional spaces with $m\ge 1$ fixed and
let $\widehat F$ be their ultraproduct. It is well known that
$$\Pi B(E_i,F_i)/{\cal U} = B(\widehat E,\widehat F)\leqno (30)$$
isometrically. Now observe that for any integer $k$ we have (by [Sm]) for
any $u\colon \ \widehat E\to M_k$, associated to a family $(u_i)_{i\in I}$
with $u_i \in B(E_i,M_k)$
$$\|u\|_{cb} = \|I_{M_k}\otimes u \|_{M_k(\widehat E)\to M_k(M_k)}.$$
By (30) it follows that
$$\|u\|_{cb} = \lim_{\cal U} \|I_{M_k}\otimes u_i\|_{M_k(E_i) \to
M_k(M_k)}.$$
Hence
$$\|u\|_{cb(\widehat E,M_k)} = \lim_{\cal U} \|u_i\|_{cb(E_i,M_k)}.$$
Equivalently
$$\eqalign{M_k((\widehat E)^*) &= \Pi M_k(E^*_i)/{\cal U}\cr
&= M_k(\Pi E^*_i/{\cal U})}$$
and we conclude that $(\widehat E)^*$ and $\Pi E^*_i/{\cal U}$ are
completely isometric.\qed

\proclaim Corollary 15. Consider $E$ in $OS_n$. The following are
equivalent:
\medskip
\item{(i)} $d_{SK}(E) = d_{SK}(E^*)=1$.
\item{(ii)} For any sequence $E_i$ in $OS_n$ such that $\widehat E = \Pi
E_i/{\cal U}$ is completely isometric to $E$ we have
$$\lim_{\cal U} d_{cb}(E,E_i) =1.$$
\item{(iii)} Same as (ii) with each $E_i$ isometric to $E$.\medskip

\pf Assume (ii) (resp. (iii)). Then the proof of Proposition~12 (resp.
Corollary~13) shows that necessarily $d_{SK}(E)=1$. By Lemma~14 it is clear
that (ii) and (iii) are self dual properties, hence we must also have
$d_{SK}(E^*)=1$. Conversely assume (i). We will use Proposition~6 together
with the identity $E\otimes_{\rm min} F = cb(F^*,E)$ valid when $F$ is
finite dimensional. By Proposition~6 and Lemma 14, if $d_{SK}(E) = 1$ we have
$$\Pi cb(E_i,E)/{\cal U} = cb(\widehat E,E)\leqno (31)$$
whenever $\widehat E$ is an ultraproduct of $n$ dimensional spaces. Now if
$d_{SK}(E^*)=1$ this implies
$$\Pi cb(E^*_i,E^*)/{\cal U} = cb(\widehat E^*, E^*)$$
hence after transposition
$$\Pi cb(E,E_i)/{\cal U} = cb(E,\widehat E).\leqno (32)$$
It is then easy to conclude:\ let $u\colon \ \widehat E\to E$ be a
complete isometry. Let $u_i\colon \ E_i\to E$ be associated to $u$ via (31)
and let $v_i\colon \ E\to E_i$ be associated to $u^{-1}$ via (32) in such a
way that $\|u\|_{cb} = \lim\limits_{\cal U} \|u_i\|_{cb}=1$ and
$\|u^{-1}\|_{cb} = \lim\limits_{\cal U} \|v_i\|_{cb}=1$.

\n Then $I_E = \lim\limits_{\cal U} u_iv_i$ hence 
we have $\lim\limits_{\cal
U} \|I_E - u_iv_i\|=0$, hence $(u_iv_i)^{-1}$ exists for $i$ large and its
norm tends to 1, so that $u^{-1}_i$ exists and (since $\lim\limits_{\cal U}\|v_i\|<\infty$) $\lim\limits_{\cal U}
\|u^{-1}_i\|<\infty$, whence $\lim\limits_{\cal U}\|u^{-1}_i-v_i\|=0$. Since all
norms are equivalent on a finite dimensional space, 
we also have
$\lim\limits_{\cal
U} \|I_E - u_iv_i\|_{cb}=0$, finally (repeating the argument with
the $cb$-norms) we obtain  $\lim\limits_{\cal U}\|u^{-1}_i-v_i\|_{cb}=0$
and we conclude that
$$\lim_{\cal U} d_{cb}(E,E_i) \le \lim_{\cal U} \|u_i\|_{cb}
\|u^{-1}_i\|_{cb} \le 1.$$
This shows that (i) $\Rightarrow$ (ii) and (ii) $\Rightarrow$ (iii) is
trivial.\qed

\n {\bf Remark.}  The 
 row and column Hilbert spaces 
$R_n$ and $C_n$ obviously satisfy the properties in Corollary 15. Also, the two dimensional spaces $\ell_1^2$
and $\ell_\infty^2$, which (see [Pa2]) admit only one operator space
structure (so that any space isometric to either one is automatically 
completely isometric to it), must clearly satisfy these properties.
At the time of this writing, 
these  are the only examples I know of spaces
 satisfying the properties in Corollary 15. 

It is natural to describe the spaces appearing in Corollary~15 as points of
continuity with respect to a weaker topology (on the
metric space $OS_n$) which can be defined as follows. For
any $k\ge 1$ and any linear map $u\colon \ E\to F$ between
operator spaces, we denote as usual $$\|u\|_k =
\|I_{M_k}\otimes E\|_{M_k(E)\to M_k(F)}.$$ Then for any
$E,F$ in $OS_n$ we define $$d_k(E,F) = \inf\{\|u\|_k
\|u^{-1}\|_k\}$$ where the infimum runs over all
isomorphisms $u\colon\ E\to  F$.

We will say that a sequence $\{E_i\}$ in $OS_n$ tends weakly to $E$ if, for
each $k\ge 1$, $\hbox{Log } d_k(E_i,E)\to 0$ when $i\to \infty$. This
notion of limit clearly corresponds to a topology (namely to the topology
associated to the metric $\tilde\delta = \sum\limits_{k\ge 1} 2^{-k} \hbox{Log }d_k$)
which we will call the weak topology. Let us say that $E,F$ are
$k$-isometric
if there is an isomorphism  $u\colon\  E\to F$ such that $I_{M_k}\otimes u$
is an isometry. Clearly this holds (by a compactness argument) iff
$d_k(E,F) =1$. Moreover (again by a compactness argument) $E$ and $F$ are
completely isometric iff they are $k$-isometric for all $k\ge 1$. This
shows that the weak topology on $OS_n$ is Hausdorff. We observe

\proclaim Proposition 16. Let $E$ and $E_i$ ($i=1,2\ldots$) be operator
spaces in $OS_n$. Then $E_i$ tends weakly to $E$ iff for any nontrivial
ultrafilter ${\cal U}$ on $\bf N$ the ultraproduct $\widehat E = \Pi
E_i/{\cal U}$ is completely isometric to $E$.

\pf Clearly if $E_i$ tends weakly to $E$ then $\widehat E$ is $k$-isometric
to $E$ for each $k\ge 1$, hence $\widehat E$ is completely isometric to
$E$. Conversely, if $E_i$ does not tend weakly to $E$, then for some $k\ge
1$ and some $\vp>0$ there is a subsequence $E_{n_i}$ such that
$$d_k(E_{n_i},E) >1+\vp \quad \hbox{for all}\quad i\ne
j.\leqno (33)$$
Let ${\cal U}$ be an ultrafilter refining this subsequence, let $\widehat
E$ be the corresponding ultraproduct, and let $u\colon \ \widehat E\to E$
be any isomorphism. Clearly there are isomorphisms $u_i\colon \ E_i\to E$
which correspond to $u$ and we have for each $k\ge 1$ (by compactness)
$\|u\|_k = \lim\limits_{\cal U} \|u_i\|_k$ and $\|u^{-1}\|_k =
\lim\limits_{\cal U} \|u^{-1}_i\|_k$. Therefore we obtain by (33)
$d_k(\widehat E,E)\ge 1+\vp$ and we conclude that
$E$ and $\widehat E$ are not completely isometric. \qed

We can now reformulate Corollary 15 as follows

\proclaim Corollary 17. Let $i\colon \ OS_n \to OS_n$ be the identity
considered as a map from $OS_n$ equipped with the weak topology to $OS_n$
equipped with the metric $d_{cb}$. Then an element $E$ in $OS_n$ is a point
of continuity of $i$ iff $d_{SK}(E) = d_{SK}(E^*)=1$.

\vskip24pt

\magnification\magstep1
\baselineskip = 18pt
\def\n{\noindent}

\overfullrule = 0pt
\def\pf{\medskip{\noindent{\bf Proof. }}}
\def\vp{\varepsilon}

\n {\bf \S 5. On the dimension of the containing matrix space.}

It is natural to try to connect our study of $d_{SK}(E)$ with a result of
Roger Smith [Sm]. Smith's result implies that for any subspace $F\subset
M_k$ we have for any operator space $X$ and for any linear map $u\colon \
X\to F$
$$\|u\|_{cb}\le \|u\|_k.\leqno (34)$$
More generally, if $T$ is any compact set and $F\subset
C(T)\otimes_{\rm min} M_k$ then we also have (34).

Now let $E$ be any finite dimensional operator space. For each $k\ge 1$ we
introduce
$$\delta_k(E) = \inf\{d_{cb}(E,F)\mid F\subset
C(T)\otimes_{\rm min} M_k\}$$
 where the infimum runs over all possible compact sets
$T$.

\n Note that
$$d_{SK}(E) = \inf_{k\ge 1} \delta_k(E).$$
Clearly if $C = \delta_k(E)$, then by (34) we have for all $u\colon \ X\to
E$
$$\|u\|_{cb} \le C\|u\|_k.\leqno (35)$$
It turns out that the converse is true:\ if (35) holds for all $X$ and all
$u\colon \ X\to E$ then necessarily $\delta_k(E)\le C$. This is contained
in the next statement which can be proved following the framework of
[P2], but  using an idea of Marius Junge [J].

\proclaim Theorem 18. Let $E$ be any operator space and let $C\ge 1$ be a
constant. Fix an integer $k\ge 1$. Then the following are
equivalent:\medskip
\item{(i)} There is a compact set $T$ and a subspace
$F\subset C(T)\otimes_{\rm min} M_k$ such that
$d_{cb}(E,F) \le C$. \item{(ii)} For any operator space
$X$ and any $u\colon \ X\to E$ we have $$\|u\|_{cb} \le
C\|u\|_k.$$ \item{(iii)} For all finite dimensional
operator spaces $X$, the same as (ii) holds.\medskip

\pf (Sketch) (i) $\Rightarrow$ (ii) is Smith's result [Sm]
as explained above. (ii)~$\Rightarrow$~(iii) is trivial.
Let us prove (iii)~$\Rightarrow$~(i). Assume (iii). Note
that (for $k$ fixed) the class of spaces of the form $C(T)
\otimes_{\rm min} M_k$ is stable by ultraproduct, since
the class of commutative unital $C^*$-algebras is stable
by ultraproducts. In particular we may and will assume
(for simplicity) that $E$ is finite dimensional. Let $G$ be
an other operator space and consider a linear map $v\colon\
E\to G$. We introduce the number $\alpha_k(v)$ as folows.
We consider all factorizations of $v$ of the form
$$E\ {\buildrel a\over
\longrightarrow} \ \ell^N_\infty \otimes_{\rm min} M_k\
{\buildrel b\over \longrightarrow}\ G$$ 
where $N\ge 1$ is an arbitrary integer, and and we set
$$\alpha_k(v) =
\inf\{\|a\|_{cb}\|b\|_{cb}\}$$ where the infimum runs over
all possible $N$ and all possible such factorizations.

 \n Now,
using Lemma~2 and an ultraproduct argument it suffices to
prove that for any $n$, any $\vp>0$ and any map $v\colon\ 
E\to M_n$, there is an integer $N\ge 1$ and a
factorization of $v$ of the form $$E\ {\buildrel a\over
\longrightarrow} \ \ell^N_\infty \otimes_{\rm min} M_k\
{\buildrel b\over \longrightarrow}\ M_n$$
with $\|a\|_{cb}\|b\|_{cb}\le C(1+\vp)\|v\|_{cb}$.

\n In other words, to conclude the proof, it suffices to
show that if (iii) holds we have for all $n$ and all $v\colon \ E\to
M_n$ $$\alpha_k(v) \le C\|v\|_{cb}.\leqno (36)$$
Now we observe that (for any $G$)  $\alpha_k$ is a norm on
$cb(E,G)$ (left to the reader, this is where the presence
of $\ell^N_\infty$ is used). Hence (36) is equivalent to a
statement on the dual norms. More precisely, (36) is
equivalent to the fact that for any $T\colon \ M_n\to E$
we have $$|tr(vT)|\le C\|v\|_{cb} \alpha^*_k(T),\leqno
(37)$$ where $\alpha^*_k$ is the dual norm to $\alpha_k$,
i.e. $$\alpha^*_k(T) = \sup\{|tr(vT)|\mid v\colon \ E\to
M_n,\  \alpha_k(v)\le 1\}.$$

\n We will now use the operator space version of the
absolutely summing norm which was first introduced in 
[ER6]. In the broader framework of [P2], these operators
are called completely 1-summing and the corresponding norm
is denoted by $\pi^o_1$. We will use a version of the
``Pietsch factorization'' for these maps which is
presented in [P2]. 
As observed by M. Junge, in the present situation 
the proof of Theorem~2.1
and Remark~2.7 in [P2] yield a factorization of $T$ 
of the following form
$$M_n\ {\buildrel w\over
\longrightarrow} \ X\
{\buildrel u\over \longrightarrow}\ E$$
where $X$ is an $n^2$-dimensional operator space and where
the maps $w$ and $u$ satisfy
$$\pi^o_1(w) \le 1\quad \hbox{and}\quad \|u\|_k \le \alpha^*_k(T).$$
Hence if (iii) holds we find
$$\pi^o_1(T) = \pi^o_1(uw) \le \|u\|_{cb} \pi^o_1(w) \le \|u\|_{cb}\le C\|u\|_k \le
C\alpha^*_k(T).$$
But since $T$ is defined on $M_n$, we clearly have by definition of
$\pi^o_1$
 $$|tr(vT)|\le \pi^o_1(vT) \le \pi^o_1(T) \|v\|_{cb},$$
hence we obtain (37).\qed

Junge's idea can also be used to obtain many  more variants. For instance
let $k\ge 1$ be fixed. Consider the following property of an operator space
$E$: There is a constant $C$ such that for any $n$ and any bounded operator $v\colon \ M_n\to M_n$ we have
$$\|v\otimes I_E\|_{M_n(E)\to M_n(E)} \le
C\|v\|_k=C\|v\|_{M_n(M_k)\to M_n(M_k)}.$$ Then this property
holds iff there is an operator space $F$ completely
isomorphic to  $E$ with $d_{cb}(E,F)\le C$ such  that for
some compact set $T$, $F$ is a quotient of a subspace of
$C(T) \otimes_{\rm min} M_k$. In particular if $k=1$ this
result answers a question that I had raised in the problem
book of the Durham meeting in July~92. Here is a sketch of
a proof: By ultraproduct arguments, an operator space $E$ 
satisfies $d_{cb}(E,F)\le C$ for some $F$  subspace of a
quotient of  $C(T) \otimes_{\rm
min} M_k$ for some $T$ iff for any integer $n$ and any maps
$v_1 \colon\  M_n^* \to E$ and $v_2 \colon\  E\to M_n$ we have
$$\alpha_k(v_2 v_1) \le C \|v_1\|_{cb} \|v_2\|_{cb}.$$
Equivalently, this holds iff
$$|tr(v_2 v_1 T) | \le C \|v_1\|_{cb} \|v_2\|_{cb}
\alpha^*_k(T) $$
for any $T \colon\  M_n \to M_n^*$. The proof can then be
completed using the factorization property of $\alpha^*_k$
as above. (See [Her] for similar
 results in the category of Banach spaces.)

\n We refer the reader to
M.~Junge's forthcoming paper for more results of this kind.

\magnification\magstep1
 \baselineskip = 18pt
 \def\n{\noindent}

 \overfullrule = 0pt
 \def\qed{{\hfill{\vrule height7pt width7pt
depth0pt}\par\bigskip}} 

\vskip24pt
 \centerline{\bf References}

\item{[AO]}   C. Akeman and P. Ostrand.
 Computing norms in group $C^*$-algebras. 
Amer. J. Math. 98 (1976), 1015-1047.

 \item{[B1]} D. Blecher. Tensor products of operator spaces II. (Preprint)
 1990.  Canadian J. Math. 44 (1992) 75-90.
 
 \item{[B2]}  $\underline{\hskip1.5in}$. The standard dual of an operator space.
  Pacific J. Math. 153 (1992) 15-30.

 \item{[BP]} D. Blecher and V. Paulsen. Tensor products of operator spaces.
 J. Funct. Anal. 99 (1991) 262-292.

\item{ [DCH]} J. de Canni\`ere   and  U. Haagerup. 
Multipliers of the Fourier algebras of some simple Lie
groups and their discrete subgroups.  Amer. J. Math.   107
(1985), 455-500.

 \item{[EH]} E. Effros and U. Haagerup. Lifting problems
and local reflexivity for $C^*$-algebras. Duke Math. J.
52 (1985) 103-128.

 \item{[EKR]} E. Effros J. Kraus and Z.J. Ruan. On two
quantized tensor products. Preprint 1992. To appear.

 \item{[ER1]}  E. Effros  and Z.J. Ruan. A new 
approach to operator spaces.
 Canadian Math. Bull.
34 (1991) 329-337.

 \item{[ER2]} $\underline{\hskip1.5in}$. On the abstract 
characterization of
 operator spaces. Proc. A.M.S. To appear.
 
 \item{[ER3]} $\underline{\hskip1.5in}$. Self duality for the Haagerup
 tensor product and Hilbert space factorization.  J. Funct. Anal. 100
(1991) 257-284.
 
 \item{[ER4]} $\underline{\hskip1.5in}$. Recent development in operator
 spaces. Current topics in Operator algebras (edited by Araki, Choda, Nakagami, Sait\^o
and Tomiyama). World Scientific, Singapore, 1991 p. 146-164.

 \item{[ER5]} $\underline{\hskip1.5in}$. Mapping spaces
and liftings for operator spaces. (Preprint) Proc. London
Math. Soc. To appear.

\item{[ER6]}  $\underline{\hskip1.5in}$. The
Grothendieck-Pietsch and Dvoretzky-Rogers Theorems for
operator spaces. (Preprint 1991) J. Funct. Anal. To
appear.

\item{[ER7]} $\underline{\hskip1.5in}$. On approximation
properties for operator spaces, International J. Math. 1
(1990) 163-187.

 \item{[H]} U. Haagerup. Injectivity and decomposition of completely
 bounded maps in ``Operator algebras and their connection with Topology and
 Ergodic Theory''. Springer Lecture Notes in Math. 1132
(1985) 170-222.

 \item{[HP]}  U. Haagerup and G. Pisier.  Bounded linear
operators between
 $C^*$-algebras. Duke Math. J. 71 (1993) To appear.
  
\item{[Hei]} S. Heinrich. Ultraproducts in Banach space
theory. J. f\"ur die reine und Angew. Math. 313 (1980)
72-104.

\item{[Her]} R. Hernandez. Espaces $L^p$, factorisation et produits
tensoriels dans les espaces de Banach.	Comptes Rendus Acad. Sci. Paris 
S\'erie A. {296} (1983) 385--388.

 \item{[J]} M. Junge. Oral communication.

 \item{[Ki]} E. Kirchberg. On subalgebras of the
CAR-algebra. (Preprint, Heidelberg, 1990) To appear.

 \item{[Kr]} J. Kraus. The slice map problem and
approximation properties. J. Funct. Anal. 102 (1991)
116-155.

\item{[LR]} J. Lindenstrauss and H. Rosenthal. 
The ${\cal L}_p$ spaces. Israel
J. Math. 7 (1969) 325-349.

 \item{[Pa1]} V. Paulsen.  Completely
bounded maps and dilations. Pitman Research Notes 146.
Pitman Longman (Wiley) 1986.

\item{[Pa2]} $\underline{\hskip1.5in}$. Representation of
Function algebras, Abstract
 operator spaces and Banach space Geometry. J.
Funct. Anal. 109 (1992) 113-129.

 \item{[P1]} G. Pisier. The operator Hilbert space $OH$,
complex interpolation and tensor norms. To appear.

 \item{[P2]}  $\underline{\hskip1.5in}$. Non-commutative vector valued $L_p$-spaces
and completely $p$-summing maps. To appear.

\item{[Ru]} Z.J. Ruan. Subspaces of $C^*$-algebras. J. Funct. Anal. 76
 (1988) 217-230.
 
\item{[Sm]} R.R. Smith. Completely bounded maps
between $C^*$-algebras. J. London Math. Soc. 27
(1983) 157-166.
 
\item{[Sz]} A. Szankowski. $B(H)$ does not have the approximation property.
Acta Math. 147 (1981) 89-108.

\item{[Ta]} M. Takesaki. Theory of Operator Algebras I.
Springer-Verlag New-York 1979.

\item{[W1]} S.Wasserman. On tensor products of certain
group $C^*$-algebras. J. Funct. Anal. 23 (1976) 239-254.

\item{[W2]}$\underline{\hskip1.5in}$.The slice map problem
 for $C^*$-algebras.
Proc. London Math. Soc. 32 (1976) 537-559.

\vskip12pt

\vskip12pt

Texas A. and M. University, College Station, TX 77843, U. S. A.

and

Universit\'e Paris 6, Equipe d'Analyse, Bo\^\i te 186,
75252 Paris Cedex 05, France

 \end

 \n {\bf Remark.} In conclusion, we would like to mention
a possible direction of research in analogy with the
"Local" Theory of Banach spaces. Consider an
$n$-dimensional operator space $E$ such that
$d_{SK}(E)=1$. The problem is to find estimates for the
number
$$N(E)= \inf\{k| \delta_k(E)\le 2\}.$$
(Here we choose $2$ for simplicity, one could put instead
$1+\vp$, for some fixed $\vp>0$.)
More explicitly, 

is there an upper bound of $N(E)$ in
terms

$$\matrix{&X\cr&\cr &w\nearrow\quad \quad \searrow u\cr&&\cr
&M_n \quad {}_{\displaystyle {\buildrel {\hbox to
30pt{\rightarrowfill}} \over T}} \quad E\cr}$$
 
 \item{[P2]} $\underline{\hskip1.5in}$. Factorizations of operator valued
 analytic functions and complex interpolation. Feitschrift in honor of
 I.~Piatetski-Shapiro (edited by Gelbart, Howe and Sarnak)
Part II,
 pp.~197-220, Weizmann Science Press of Israel, IMCP, 1990.

 \item {[P6]} $\underline{\hskip1.5in}$. Espace de Hilbert
d'op\'erateurs et interpolation complexe. Comptes Rendus
Acad. Sci. Paris S\'erie I, 316 (1993)  47-52.

\item {[P7]} $\underline{\hskip1.5in}$. Sur les
oprateurs factorisables par $OH$.  Comptes Rendus Acad.
Sci. Paris S\'erie I, 316 (1993) 165-170.

\item{[P3]} $\underline{\hskip1.5in}$. Remarks on
complemented subspaces of $C^*$-algebras. Proc. Roy. Soc.
Edinburgh, 121A (1992) 1-4.

\item{[Be]} J. Bergh. On the relation between the two
complex methods of interpolation. Indiana Univ. Math.
Journal 28 (1979) 775-777.

 \item{[BSp]} M. Bo$\dot{z}$ejko and R. Speicher.
An example of a generalized Brownian motion.
Commun. Math. Phys. 137 (1991) 519-531. 

 \item{[Ca]} A. Calder\'on. Intermediate spaces and interpolation, the
 complex method. Studia Math. 24 (1964) 113-190.

\item{[Co]} A. Connes. Classification of injective
factors, Cases
$II_1,II_\infty,III_\lambda,\lambda\neq 1$. Ann.
Math. 104 (1976) 73-116.

\item{ [DCH]} J. de Canni\`ere   and  U. Haagerup. 
Multipliers of the Fourier algebras of some simple Lie
groups and their discrete subgroups.  Amer. J. Math.   107
(1985), 455-500.
 \item{[Gl1]} E.Gluskin. The diameter of the Minkowski
compactum is roughly equal to $n$. Funct. Anal. Appl. 15
(1981) 72-73.

\item{[Gl2]} $\underline{\hskip1.5in}$. Probability in
the geometry of Banach spaces. Proceedings I.C.M.
Berkeley, U.S.A.,1986 Vol 2, 924-938 (Russian),
Translated in "Ten papers at the I.C.M. Berkeley 1986"
Translations of the A.M.S.(1990)

\item {[GK]} J. Grosberg and M. Krein.  Sur la
d\'ecomposition des fonctionnelles en composantes
positives.  Doklady Akad. Nauk SSSR {\bf 25} (1939),
723-726.

\item{[HP1]} U. Haagerup and G. Pisier. Factorization of analytic functions
 with values in non-commutative $L^1$-spaces and applications. Canadian J.
 Math. 41 (1989) 882-906.
 
 \item{[HP2]} $\underline{\hskip1.5in}$.  Bounded linear
operators between
 $C^*$-algebras. To appear.
 
\item{[HeP]} A.Hess and G.Pisier. On the $K$-functional
for the couple of operator spaces  $(B(c_0),B(\ell_1))$.
In preparation.

 \item{[K]} O. Kouba. Interpolation of injective or projective tensor
 products of Banach spaces. J. Funct. Anal. 96 (1991) 38-61.
`

\item{[T]} J.Tomiyama. On the projection of norm one in
$W^*$-algebras. Proc. Japan Acad. 33 (1957) 608-612.

\item {[P10]} $\underline{\hskip1.5in}$. Noncommutative
 vector valued integration and operator
$p$-summing maps. In preparation.

\item{[LT1]} J.Lindenstrauss and L.Tzafriri. Cla
\item{[L]} C. Le Merdy. Analytic factorizations and completely bounded
 maps. Preprint. June 1992.
 ssical
Banach spaces, vol. I. Sequence spaces. Springer Verlag,
Berlin 1976.

\item{[LT2]} $\underline{\hskip1.5in}$. Classical
Banach spaces, vol. II. Function spaces. Springer Verlag,
Berlin 1979.

\item{[LPP]}  F. Lust-Piquard  and G. Pisier.   Non
commutative Khintchine and Paley inequalities.  Arkiv fr
Mat.  29 (1991) 241-260.
\item {[P5]} $\underline{\hskip1.5in}$. The volume of
convex bodies and Banach space geometry.
Cambridge Univ. Press, 1989.

\bye